\documentclass[10pt]{article}
\usepackage{graphicx}
\usepackage{amsthm}
\usepackage{amsmath,amssymb,wasysym}
\usepackage[left=3cm,right=4cm, top=3.3cm, bottom=3cm]{geometry}
\newtheorem{theo}{Theorem}[section]

\newtheorem{coro}[theo]{Corolary}
\newtheorem{prop}[theo]{Proposition}

\newtheorem{remark}[theo]{Remark}
\def\proof {{\noindent \bf{Proof:\hspace{4pt}}}}
\def\endproof{\hfill$\square$\vspace{6pt}}
\numberwithin{equation}{section}
\title{
{\bf\Large The Symmetric Regularized-Long-Wave Equation: Ill-posedness and Long Period Limit}}
\author{
{\bf\large Carlos Banquet}\footnote{Email:  cabanquet@hotmail.com, cbanquet@correo.unicordoba.edu.co}\hspace{2mm}
{\bf\large}\vspace{1mm}\\
{\it\small Departamento de Matem\'aticas y Estad\'istica}\\
{\it\small  Universidad de C\'ordoba}\\
{\it\small Carrera 6 No. 76-103, Monter\'ia, C\'ordoba, Colombia}\vspace{3mm}}
\date{}
\begin{document}
\maketitle
\begin{abstract}
In the present work we obtain two important results for the Symmetric Regulraized-Long-Wave equation. First we prove that the initial value problem for this equation  is ill-posed for data in  $H^s(\mathbb{R})\times H^{s-1}(\mathbb{R}),$  if $s< 0,$ in the sense that the flow-map cannot be continuous at the origin from $H^s(\mathbb{R})\times H^{s-1}(\mathbb{R})$ to even $(\mathcal{D}'(\mathbb{R}))^2.$ We also establish an exact theory of convergence of the periodic solutions  to the continuous one, in Sobolev spaces, as the period goes to infinity.
\end{abstract}


%%%%%%%%%%%%% INTRODUCTION %%%%%%%%%%%%%%% INTRODUCTION %%%%%%%%%%%%%%%%%%%%%%%%%%%%%%%%%%%%%%%%
\section{Introduction}
In this work we study the Symmetric  Regularized-Long-Wave  equation (SRLW equation henceforth)
\begin{equation} \label{ecuabasica}
u_{tt}+u_{xx}+(uu_x)_t-u_{xxtt}=0,
\end{equation}
where $u$ is a real-valued function and the subscripts  denote the derivative with respect to the  spatial variable $x$ and time  $t.$ The Symmetric Regularized-Long-Wave equation is a model for the weakly nonlinear ion acoustic and space-charge waves. This equation was introduced by Seyler and Fenstermacher in \cite{SeylerFen}, where a  weakly nonlinear analysis of the cold-electron fluid equation is made. 
The equation $(\ref{ecuabasica})$ has the equivalent form 
\begin{equation} \label{ecuaequiv}
\left \{ 
\begin{aligned}
&u_t-u_{xxt}+uu_x-v_x=0,\\
&v_t-u_x=0,\\
\end{aligned} \right.
\end{equation}
for all  $t>0$ and $x\in \mathbb{R}.$ This equivalent system has four conservation laws
\begin{equation*}
\begin{split}
E(u,v)&= \int (uv -\frac{1}{6} u^3)dx,\ \ \ V(u,v)= \frac{1}{2}\int (u^2+u_x^2+v^2)dx,\\
I_1(u,v)&=\int u\ dx\ \ \  \text{and}\ \ \ I_2(u,v)=\int v\ dx.
\end{split}
\end{equation*}
The SRLW equation have been studied for different authors from various point of view, for example Chen in \cite{ChenLin} established the stability and instability of the solitary wave solutions associated to the generalized SRLW equation using the framework given by Grillakis, Shatah and Strauss \cite{grillakis2, grillakis3}.  See also Chen and Li \cite{ChenLi}, Fei \cite{Feixu}, and Zhang \cite{Zhang}, where different technics were used to obtain explicit solitary and periodic solutions for  generalized versions of the SRLW equation.\\

For later use, let us define $D_x=\frac{1}{i}\partial_x.$ Then (\ref{ecuaequiv}) can be written as
\begin{equation} \label{ecuawithi}
\left \{ 
\begin{aligned}
&iu_t=\varphi(D_x)\left(\tfrac{1}{2}u^2-v\right),\\
&iv_t=-\psi(D_x)u,\\
&u(x,0)=u_0(x), \ \ \ v(x,0)=v_0(x),
\end{aligned} \right.
\end{equation}
where $\varphi(D_x)$ and $\psi(D_x)$ are  the Fourier multiplier operators define by
\[
\widehat{\varphi(D_x)u}(\xi)=\varphi(\xi)\widehat{u}(\xi), \ \ \ \widehat{\psi(D_x)u}(\xi)=\psi(\xi)\widehat{u}(\xi), \ \ \text{with}\  \ \ \varphi(\xi)=\frac{\xi}{1+|\xi|^2}, \ \ \psi(\xi)=\xi.
\]
Quite recently, Banquet in \cite{Banquet1} prove that the initial value problem (\ref{ecuawithi}) is globally well-posed in the Sobolev spaces $H^s(\mathbb{R})\times H^{s-1}(\mathbb{R})$ for $s\geq 0.$ More precisely,
\begin{theo}\label{TheoBanquet}
Suppose $s\geq 0,$ then for all $(u_0,v_0)\in H^s(\mathbb{R})\times H(\mathbb{R})^{s-1}$ there exists a positive $T=T(\|u_0\|_{H^s},\|v_0\|_{H^{s-1}})$  and a unique solution of (\ref{ecuawithi}) on the interval $[-T,T]$, such that $(u,v) \in C ([-T,T];H^s(\mathbb{R})\times H^{s-1}(\mathbb{R})).$\\
Furthermore, for  $R>0,$ let $\mathcal{B}_R$ denote the ball of radius $R$ centered at the origin in $H^s(\mathbb{R})\times H^{s-1}(\mathbb{R})$ and let $T=T(R)>0$ denote a uniform existence time for the initial value problem (\ref{ecuawithi}) with $(u_0,v_0)\in \mathcal{B}_R.$ Then the correspondence $(u_0,v_0)\mapsto (u(t),v(t))$ that associates to $(u_0,v_0)$ the solution $(u(t),v(t))$ of the IVP (\ref{ecuawithi}) with initial data $(u_0,v_0)$ is a real analytic mapping from $\mathcal{B}_R$ to $C ([-T,T];H^s(\mathbb{R})\times H^{s-1}(\mathbb{R})).$\\
\end{theo}
To obtain this result Banquet performed a Picard iteration on an adequate space and took advantage of Lemma 1  established by Bona and Tzvetkov in \cite{bonaTzvetkov}. On the periodic case a similar result as in Theorem \ref{TheoBanquet} was obtained in Banquet \cite{Banquet1}, in this case the proof in based on Lemma 3.1 established by Roum\'egoux in \cite{Roume1}. The theory given by Banquet improves the earlier work of Chen \cite{ChenLin}, where a global well-posedness result was proved in $H^1(\mathbb{R})\times L^2(\mathbb{R}).$\\

Recently, motivated by the idea introduced by Bejenaru and Tao \cite{BejTao1}, Molinet and Vento \cite{MolVen1} obtained a sharp ill-posedness result for the KdV-Burgers equation in $H^s(\mathbb{R}),$ $s<-1$ in the sense that the flow map $u_0\mapsto u(t)$ cannot be continuous from $H^s(\mathbb{R}),$  to even $\mathcal{D}'(\mathbb{R})$ at any fixed $t>0$ small enough. It is worth to note that Panthee in \cite{Panthee1} based on the ideas on Molinet and Vento proved that the IVP for the BBM equation is ill-posed for data in $H^s(\mathbb{R})$ if $s<0,$ in the sense described above.  Following the scheme presented by Molinet and Vento, in this paper we obtain the next result.
\begin{theo}\label{TheoBanquet2}
Let $s<0,$ then the initial value problem (\ref{ecuawithi}) is ill-posed in $H_{per}^s(\mathbb{R})\times H_{per}^{s-1}(\mathbb{R})$ in the following sense: There exists $T>0$ such that for any $0<t<T,$ the flow map $(u_0,v_0)\mapsto(u(t),v(t))$ constructed in Theorem \ref{TheoBanquet}  is discontinuous at the origin from $L^2(\mathbb{R})\times H^{-1}(\mathbb{R})$ endowed with the  topology induced by $H_{per}^s(\mathbb{R})\times H_{per}^{s-1}(\mathbb{R})$ into $\mathcal{D}'(\mathbb{R})\times \mathcal{D}'(\mathbb{R}).$
\end{theo}
We can also prove a very similar result as in  Theorem \ref{TheoBanquet2} for the SRLW equation in the periodic setting (See Remark \ref{remark1} below). At this point we want to refer to the work of Molinet and Vento \cite{MolVen2}, where the sharp ill-posedness result for the KdV-Burgers equation is proved in the periodic case. See also Panthee \cite{Panthee1}, where a similar result is obtained for the periodic  BBM equation.\\

The second part of the paper is dedicated to study another important property of the regularized equations, namely, the convergence of the periodic solutions to the continuous one on Sobolev spaces. Here we propose a scheme for approximating a localized disturbance by a spatially periodic evolution. More precisely we study the SRLW equation with two initial data, one vanishing at infinity and another periodic with large period.\\

 Next, to give the study focus, the main result of second part  the paper is stated informally.
 \begin{theo}
 Let $(u,v)$ the solution of the SRLW equation (\ref{ecuaequiv}) corresponding to the initial condition $(u(x,0),v(x,0))=(\psi(x),\phi(x))$ which is sufficiently nice. Let $\mathcal{P}_l(\psi)$ and $\mathcal{P}_l(\phi)$ appropriately version of $\psi$ and $\phi,$ respectively. Consider the periodic solution of $(u_l,v_l),$ with initial data $(\mathcal{P}_l(\psi), \mathcal{P}_l(\phi)).$  Then, when both solutions are restricted to the spatial interval $(-l,l),$ their difference satisfies
 \[\lim_{l\to\infty}\|(u_l(t),v_l(t))-(u(t),v(t))\|_{X^{1,2}_l}=0,\] 
 uniformly on compact time intervals.
 \end{theo}
 This theory of convergence provides a discrete method to approximate the solutions of the SRLW equation when the initial data lies in some Sobolev space on the line. As far as we know it does not exist in the literature these kind of results about convergence and ill-posedness, we hope that these new theories would be of interest.\\ 
 
This paper is organized as follows: In Section 2 we introduced some notations to be used throughout the whole article; in Section 3, we obtain the results of ill-posedness in the continuous  and periodic setting; in Section 4, we show the convergence of the  periodic solutions to the continuous one. 
%%%%%%%%%%%%%%%%%%%%%%%%%%% NOTATION %%%%%%% NOTATION %%%%%%%%%%%%%%%%%%%%%%%%%%%%%%%%%%%%%
\section{Notation and preliminaries}
We denote by $\mathcal{F}_x(f)(\xi)$ or $\widehat{f}(\xi),$ the Fourier transform of $f$ in $x$ variable
\[\widehat{f}(\xi):=\widehat{f}(\xi):=\frac{1}{\sqrt{2\pi}}\int_{\mathbb{R}}e^{-ix\xi}f(x)dx.\]
We use $H^s(\mathbb{R})$ to denote the $L^2-$based Sobolev space of order $s$ with norm 
\[\|f\|_{H^s}=\left(\int_{\mathbb{R}}(1+|\xi|^2)^s|\widehat{f}(\xi)|^2d\xi\right)^{\frac 1{2}}.\]
Similarly, if $I$ is an interval in $\mathbb{R},$ $L^p,$ with $p\geq 1$ represents the usual Lebesgue space with the usual norm $\|\cdot\|_{L^p(I)}.$ In particular, we denote $\|\cdot \|_{L^{\infty}(\mathbb{R})}$ by $\|\cdot\|_{\infty}.$ For any nonnegative integer  $m,$ the space 
\[W^{m,2}(I)=\left\{f,f',...,f^{(m)}\in L^2(I): \int_I\left(|f(x)|^2+|f^{(m)}(x)|^2\right)dx<\infty\right\}\] 
is a Hilbert space with norm defined by 
\[\|f\|_{W^{m,2}(I)}=\left(\int_I\left(|f(x)|^2+|f^{(m)}(x)|^2\right)dx\right)^{\frac 1{2}}.\] 
If $m=0,$ we denote $W^{0,2}(I)=L^2(I).$ We also use the Hilbert space $X^{m,2}(I)=W^{m,2}(I)\times W^{m-1,2}(I)$ with the natural norm.
\[\|(f,g)\|_{X^{m,2}}=\|f\|_{W^{m,2,}}+\|g\|_{W^{m-1,2}}.\]

Now,  we introduced an analogous periodic function space 
\[H^s_l=\left\{f(x)=\sum_{n\in\mathbb{Z}}f_ne^{\frac{in\pi}{l}}: f_n=\overline{f_{-n}}\in\mathbb{C}, \ \sum_{n\in\mathbb{Z}}(1+n^2)^s|f_n|^2<\infty\right\}\]
with the norm 
\[\|f\|_{H^s_l}=\left(2l\sum_{n\in\mathbb{Z}}\left(1+\left|\frac{n\pi}{l}\right|^{2s}\right)|f_n|^2\right)^{\frac 1{2}}.\]
If $s$ is an integer, the norm $\|f\|_{H^s_l}$ has an alternative representation
\[\|f\|_{H^s_l}=\left(\int_I\left(|f(x)|^2+|f^{(s)}(x)|^2\right)dx\right)^{\frac 1{2}},\]
which coincides with the norm $\|f\|_{W^{s,2}(I)}$ when $f$ is restricted to the interval $(-l,l).$  When $I=(-l,l),$ sometimes we write, $L^p_l,$  $X^{m,2}_l$ and $W^{m,2}_l,$ instead of $L^p(-l,l),$  $X^{m,2}(-l,l)$ and $W^{m,2}(-l,l),$ respectively.\\  

Various constants whose exact value are immaterial will be denote by $C.$ We also use the notation $A\apprle B$ (respectively, $A\apprge B$) if  there exist a positive constant  $C$ such that $A\leq CB$ (respectively, $A\geq CB$) and   $A\sim B$ means that $A\apprle B\apprle A.$\\  
%%%%%%%%%%%%%%% Ill-posedness %%%%%%%%%%%%%% Ill- posedness %%%%%%%%%%%%%%%%%%%%%%%%%%%%%%%%%%%%%%%%%
\section{Ill-posedness result}
In this section we obtain an ill-posedness result for the system (\ref{ecuaequiv}) with initial data $(u(0),v(0))=(u_0,v_0)$ on the periodic and continuous  setting.  For this we use the ideas established by Molinet  and Vento \cite{MolVen1, MolVen2} and Panthee \cite{Panthee1}. Let us start with the continuous case. Indeed, for $(u_0,v_0)\in X^s(\mathbb{R}):=H^s(\mathbb{R})\times H^{s-1}(\mathbb{R})$ consider the Cauchy problem
\begin{equation} \label{rewriteeq}
\left \{ 
\begin{aligned}
&iu_t=\varphi(D_x)\left(\tfrac 1{2}u^2-v\right),\\
&iv_t=-\psi(D_x)u\\
&(u(x,0),v(x,0))=(\epsilon u_0(x), \epsilon v_0(x)).
\end{aligned} \right.
\end{equation}
 where $\epsilon >0$ is a parameter. The solution $\overrightarrow{u^{\epsilon}}=(u^{\epsilon},v^{\epsilon})$ of (\ref{rewriteeq}) depends on the parameter $\epsilon.$ Solving the linear problem 
 \begin{equation*}
\left \{ 
\begin{aligned}
&iu_t=-\varphi(D_x)\left(v\right),\\
&iv_t=-\psi(D_x)u,\\
&(u(0),v(0))=(\epsilon u_0,\epsilon v_0)
\end{aligned} \right.
\end{equation*}
we get the solution $(u^{\epsilon}(t), v^{\epsilon}(t))=\epsilon S(t)(u_0,v_0),$ where 
\[
\left(\begin{array}{c}\widehat{u^{\epsilon}}(t,\xi)\\ \widehat{v^{\epsilon}}(t,\xi)\end{array}\right)=\epsilon \left(
\begin{array}{ccc}
\cos(\alpha(\xi)t)&\frac{i}{\sqrt{1+|\xi|^2}}\sin(\alpha(\xi)t)\\
i\sqrt{1+|\xi|^2}\sin(\alpha(\xi)t)&\cos(\alpha(\xi)t)\\
\end{array}
\right)\left(\begin{array}{c}\widehat{u}_0(\xi)\\ \widehat{v}_0(\xi)\end{array}\right)
\]
with  $\alpha(\xi)=\frac{\xi}{\sqrt{1+|\xi|^2}}.$ Then, (\ref{rewriteeq}) may be rewritten as the integral equation 
\[\overrightarrow{u^{\epsilon}}(x,t)=\epsilon S(t)(u_0,v_0)(x)-i\int_0^tS(t-\tau)G\left[\overrightarrow{u^{\epsilon}}(x,\tau)\right]d\tau\]
where $G$ is given by $G(u,v)=\left(\varphi(D_x)\left(\frac{1}{2}u^2\right),0\right).$ Then, we get 
\[\left. \frac{\partial \overrightarrow{u^{\epsilon}}(x,t)}{\partial\epsilon}\right |_{\epsilon=0}=S(t)(u_0,v_0)(x)\]
and
\[\left. \frac{\partial^2 \overrightarrow{u^{\epsilon}}(x,t)}{\partial\epsilon^2}\right |_{\epsilon=0}=-2i\int_0^tS(t-\tau)G\left[S(\tau)(u_0,v_0)(x)\right]d\tau.\]
%%%%%%%%%%%%%%% PROOF OF PRINCIPAL THEOREM OF THIS SECTION %%%%%%%%%%%%%%%%%%%%%%%%%%%%%%
\section*{Proof of Theorem \ref{TheoBanquet2}}
Let $N\gg 1$ and define $(\phi_N,\psi_N)$ via the Fourier transform as
\[\widehat{\phi_N}(\xi)=\chi_{I_N}(\xi)+\chi_{I_N}(-\xi) \ \ \text{and}\ \ \widehat{\psi_N}(\xi)=\sqrt{1+|\xi|^2}\widehat{\phi_N}(\xi),\]
where $I_N=[N-1,N+1].$ Simple calculations show that $\|\phi_N\|_{L^2(\mathbb{R})}\sim 1,$ $\|\psi_N\|_{H^{-1}(\mathbb{R})}\sim 1,$ $\|\phi_N\|_{H^{s}(\mathbb{R})}\rightarrow 0$ and $\|\psi_N\|_{H^{s-1}(\mathbb{R})}\rightarrow 0,$ if $s<0.$ 
As pointed out in the introduction of this section, the second iteration in the Picard scheme is the following,
\[
I_2((h,k),(h,k),t)=2\int_0^tS(t-\tau)G\left[S(\tau)(h,k)(x)\right]d\tau.
\]
Now, using $(\phi_N,\psi_N)$ in place of $(h,k)$ and computing the Fourier transform in $x,$ we obtain $\mathcal{F}_x\left(I_2((\phi_N,\psi_N),(\phi_N,\psi_N),t)(\xi)=\left(A(\xi,N,t), B(\xi, N, t)\right)\right),$ where
\[A(\xi,N,t)=\int_0^t\varphi(\xi)\cos(\alpha(\xi)(t-\tau))\mathcal{F}_x\left(K^2(\phi_N,\psi_N)\right)(\xi,\tau)d\tau,\]
\[B(\xi,N,t)=\int_0^ti\sqrt{1+|\xi|^2}\varphi(\xi)\sin(\alpha(\xi)(t-\tau))\mathcal{F}_x\left(K^2(\phi_N,\psi_N)\right)(\xi,\tau)d\tau\]
and
\[\widehat{K(f,g)}(\xi,t)=\cos(\alpha(\xi)t)\widehat{f}(\xi)+\frac{i}{\sqrt{1+|\xi|^2}}\sin(\alpha(\xi)t)\widehat{g}(\xi).\]
After some simple calculations we have that
\begin{align*}
A&(\xi,N,t)=\int_0^t\varphi(\xi)\cos(\alpha(\xi)(t-\tau))\int_{\mathbb{R}}\cos(\alpha(\xi-\eta)\tau)\cos(\alpha(\eta)\tau)\widehat{\phi_N}(\xi-\eta)\widehat{\phi_N}(\eta)d\eta d\tau\\
-&\int_0^t\varphi(\xi)\cos(\alpha(\xi)(t-\tau))\int_{\mathbb{R}}\frac{\sin(\alpha(\xi-\eta)\tau)\sin(\alpha(\eta)\tau)}{\sqrt{1+|\xi|^2}\sqrt{1+|\xi-\eta|^2}}\widehat{\psi_N}(\xi-\eta)\widehat{\psi_N}(\eta)d\eta d\tau\\
=&\varphi(\xi)\cos(\alpha(\xi)t)\int_{\mathbb{R}}\left[J_1(\xi,\eta,t)\widehat{\phi_N}(\xi-\eta)\widehat{\phi_N}(\eta)-\frac{J_2(\xi,\eta,t)}{\sqrt{1+|\eta|^2}\sqrt{1+|\xi-\eta|^2}}\widehat{\psi_N}(\xi-\eta)\widehat{\psi_N}(\eta)\right]d\eta\\
+&\varphi(\xi)\sin(\alpha(\xi)t)\int_{\mathbb{R}}\left[J_3(\xi,\eta,t)\widehat{\phi_N}(\xi-\eta)\widehat{\phi_N}(\eta)-\frac{J_4(\xi,\eta,t)}{\sqrt{1+|\eta|^2}\sqrt{1+|\xi-\eta|^2}}\widehat{\psi_N}(\xi-\eta)\widehat{\psi_N}(\eta)\right]d\eta,\\
\end{align*}
where
\[J_1=J_1(\xi,\eta,t)=\int_0^t\cos(\alpha(\xi)\tau)\cos(\alpha(\xi-\eta)\tau)\cos(\alpha(\eta)\tau)d\tau,\]
\[J_2=J_2(\xi,\eta,t)=\int_0^t\cos(\alpha(\xi)\tau)\sin(\alpha(\xi-\eta)\tau)\sin(\alpha(\eta)\tau)d\tau,\]
\[J_3=J_3(\xi,\eta,t)=\int_0^t\sin(\alpha(\xi)\tau)\cos(\alpha(\xi-\eta)\tau)\cos(\alpha(\eta)\tau)d\tau,\]
\[ J_4=J_4(\xi,\eta,t)=\int_0^t\sin(\alpha(\xi)\tau)\sin(\alpha(\xi-\eta)\tau)\sin(\alpha(\eta)\tau)d\tau.\]
Since $\widehat{\psi_N}(\xi)=\sqrt{1+|\xi|^2}\widehat{\phi_N}(\xi),$ we arrive at 
\begin{equation}\label{Avalue}
\begin{aligned}
A(\xi,N,t)=&\varphi(\xi)\cos(\alpha(\xi)t)\int_{\mathbb{R}}\left[J_1(\xi,\eta,t)-J_2(\xi,\eta,t)\right]\widehat{\phi_N}(\xi-\eta)\widehat{\phi_N}(\eta) d\eta\\
+&\varphi(\xi)\sin(\alpha(\xi)t)\int_{\mathbb{R}}\left[J_3(\xi,\eta,t)-J_4(\xi,\eta,t)\right]\widehat{\phi_N}(\xi-\eta)\widehat{\phi_N}(\eta) d\eta.
\end{aligned}
\end{equation}
Now, we calculate the value of $B(\xi,N, t).$
\begin{align*}
B&(\xi,N,t)=i\sqrt{1+|\xi|^2}\int_0^t\varphi(\xi)\sin(\alpha(\xi)(t-\tau))\int_{\mathbb{R}}\cos(\alpha(\xi-\eta)\tau)\cos(\alpha(\eta)\tau)\widehat{\phi_N}(\xi-\eta)\widehat{\phi_N}(\eta)d\eta d\tau\\
-&i\sqrt{1+|\xi|^2}\int_0^t\varphi(\xi)\sin(\alpha(\xi)(t-\tau))\int_{\mathbb{R}}\frac{\sin(\alpha(\xi-\eta)\tau)\sin(\alpha(\eta)\tau)}{\sqrt{1+|\xi|^2}\sqrt{1+|\xi-\eta|^2}}\widehat{\psi_N}(\xi-\eta)\widehat{\psi_N}(\eta)d\eta d\tau\\
=&i\sqrt{1+|\xi|^2}\varphi(\xi)\sin(\alpha(\xi)t)\int_{\mathbb{R}}\left[J_1(\xi,\eta,t)\widehat{\phi_N}(\xi-\eta)\widehat{\phi_N}(\eta)-\frac{J_2(\xi,\eta,t)\widehat{\psi_N}(\xi-\eta)\widehat{\psi_N}(\eta)}{\sqrt{1+|\eta|^2}\sqrt{1+|\xi-\eta|^2}}\right]d\eta\\
-&i\sqrt{1+|\xi|^2}\varphi(\xi)\cos(\alpha(\xi)t)\int_{\mathbb{R}}\left[J_3(\xi,\eta,t)\widehat{\phi_N}(\xi-\eta)\widehat{\phi_N}(\eta)-\frac{J_4(\xi,\eta,t)\widehat{\psi_N}(\xi-\eta)\widehat{\psi_N}(\eta)}{\sqrt{1+|\eta|^2}\sqrt{1+|\xi-\eta|^2}}\right]d\eta.\\
\end{align*}
Again, using the fact that $\widehat{\psi_N}(\xi)=\sqrt{1+|\xi|^2}\widehat{\phi_N}(\xi),$ we get that
\begin{equation}\label{Bvalue}
\begin{aligned}
B(\xi,N,t)=&i\sqrt{1+|\xi|^2}\varphi(\xi)\sin(\alpha(\xi)t)\int_{\mathbb{R}}\left[J_1(\xi,\eta,t)-J_2(\xi,\eta,t)\right]\widehat{\phi_N}(\xi-\eta)\widehat{\phi_N}(\eta) d\eta\\
-&i\sqrt{1+|\xi|^2}\varphi(\xi)\cos(\alpha(\xi)t)\int_{\mathbb{R}}\left[J_3(\xi,\eta,t)-J_4(\xi,\eta,t)\right]\widehat{\phi_N}(\xi-\eta)\widehat{\phi_N}(\eta) d\eta.
\end{aligned}
\end{equation}
From the definition of $\widehat{\phi},$ (\ref{Avalue}) and (\ref{Bvalue}), become
\begin{equation}\label{Avalue2}
\begin{aligned}
A(\xi,N,t)=\varphi(\xi)\cos(\alpha(\xi)t)&\int_{B_{\xi}}\left[J_1(\xi,\eta,,t)-J_2(\xi,\eta,t)\right]d\eta\\
&+\varphi(\xi)\sin(\alpha(\xi)t)\int_{B_{\xi}}\left[J_3(\xi,\eta,t)-J_4(\xi,\eta,t)\right]d\eta
\end{aligned}
\end{equation}
and
\begin{equation}\label{Bvalue2}
\begin{aligned}
B(\xi,N,t)=i\sqrt{1+|\xi|^2}\varphi(\xi)&\sin(\alpha(\xi)t)\int_{B_{\xi}}\left[J_1(\xi,\eta,t)-J_2(\xi,\eta,t)\right]d\eta\\
 &-i\sqrt{1+|\xi|^2}\varphi(\xi)\cos(\alpha(\xi)t)\int_{B_{\xi}}\left[J_3(\xi,\eta,t)-J_4(\xi,\eta,t)\right] d\eta,
\end{aligned}
\end{equation}
where
\[B_{\xi}=\{\eta: \eta\in\left[I_N\cup (-I_N)\right]\ \ \text{and}\ \ \xi- \eta\in \left[I_N\cup (-I_N)\right] \}.\]
Define
\[\theta_1(\xi,\eta)=\alpha(\xi)+\alpha(\xi-\eta)+\alpha(\eta), \ \ \theta_2(\xi,\eta)=\alpha(\xi)-\alpha(\xi-\eta)+\alpha(\eta),\]
\[\theta_3(\xi,\eta)=\alpha(\xi)-\alpha(\xi-\eta)-\alpha(\eta), \ \ \theta_4(\xi,\eta)=\alpha(\xi)+\alpha(\xi-\eta)-\alpha(\eta).\]
Then, doing some calculus we obtain that
\begin{align*}
4J_1(\xi,\eta,t)=&\frac{\sin\left[\theta_1(\xi,\eta)t\right]}{\theta_1(\xi,\eta)} +\frac{\sin\left[\theta_2(\xi,\eta)t\right]}{\theta_2(\xi,\eta)}+\frac{\sin\left[\theta_3(\xi,\eta)t\right]}{\theta_3(\xi,\eta)} +\frac{\sin\left[\theta_4(\xi,\eta)t\right]}{\theta_4(\xi,\eta)},
\end{align*}
\begin{align*}
-4J_2(\xi,\eta,t)=&\frac{\sin\left[\theta_1(\xi,\eta)t\right]}{\theta_1(\xi,\eta)} -\frac{\sin\left[\theta_2(\xi,\eta)t\right]}{\theta_2(\xi,\eta)}+\frac{\sin\left[\theta_3(\xi,\eta)t\right]}{\theta_3(\xi,\eta)} -\frac{\sin\left[\theta_4(\xi,\eta)t\right]}{\theta_4(\xi,\eta)},
\end{align*}
\begin{align*}
-4J_3(\xi,\eta,t)=&\frac{\cos\left[\theta_1(\xi,\eta)t\right]-1}{\theta_1(\xi,\eta)} +\frac{\cos\left[\theta_2(\xi,\eta)t\right]-1}{\theta_2(\xi,\eta)}+\frac{\cos\left[\theta_3(\xi,\eta)t\right]-1}{\theta_3(\xi,\eta)} +\frac{\cos\left[\theta_4(\xi,\eta)t\right]-1}{\theta_4(\xi,\eta)},
\end{align*}
\begin{align*}
4J_4(\xi,\eta,t)=&\frac{\cos\left[\theta_1(\xi,\eta)t\right]-1}{\theta_1(\xi,\eta)} -\frac{\cos\left[\theta_2(\xi,\eta)t\right]-1}{\theta_2(\xi,\eta)}+\frac{\cos\left[\theta_3(\xi,\eta)t\right]-1}{\theta_3(\xi,\eta)} -\frac{\cos\left[\theta_4(\xi,\eta)t\right]-1}{\theta_4(\xi,\eta)}.
\end{align*}
We can therefore conclude that
\begin{equation}\label{i1menosi2}
2\left[J_1-J_2\right]=\frac{\sin\left[\theta_1(\xi,\eta)t\right]}{\theta_1(\xi,\eta)} +\frac{\sin\left[\theta_3(\xi,\eta)t\right]}{\theta_3(\xi,\eta)}
\end{equation}
and 
\begin{equation}\label{i3menosi4}
-2\left[J_3-J_4\right]=\frac{\cos\left[\theta_1(\xi,\eta)t\right]-1}{\theta_1(\xi,\eta)} +\frac{\cos\left[\theta_3(\xi,\eta)t\right]-1}{\theta_3(\xi,\eta)}.
\end{equation}
Now, we move to find a lower bound for $\|I_2((\phi_N,\psi_N),(\phi_N,\psi_N),t)\|_{X^s}.$ The main contribution to this norm comes from combination of frequencies such that $|\theta_1(\xi,\eta)|$ and $|\theta_3(\xi,\eta)|$ are small. In fact, using (\ref{Avalue2}), (\ref{Bvalue2}), (\ref{i1menosi2}) and (\ref{i3menosi4}), it follows that
\begin{align}\label{IneqInt}
\|I_2((\phi,\psi)&,(\phi,\psi),t)\|^2_{X^s}=\int_{\mathbb{R}}(1+|\xi|^2)^s|A(\xi,N,t)|^2d\xi+\int_{\mathbb{R}}(1+|\xi|^2)^{s-1}|B(\xi,N,t)|^2d\xi\notag \\
&\apprge \int_{\mathbb{R}}(1+|\xi|^2)^s|\varphi(\xi)|^2\left|e^{i\alpha(\xi)t}\int_{B_{\xi}}(J_1-J_2)d\eta-ie^{i\alpha(\xi)t}\int_{B_{\xi}}(J_3-J_4)d\eta\right|^2d\xi \notag\\
&\apprge \int_{|\xi|\sim o(1)}(1+|\xi|^2)^s|\varphi(\xi)|^2\left|\int_{A_{\xi}}(J_1-J_2)d\eta\right|^2d\xi
\end{align}
where 
\[A_{\xi}=\{\eta:\eta\in I_N, \ \xi-\eta\in -I_N\ \ \text{or}\ \ \xi-\eta\in I_N, \ \ \eta\in -I_N\}.\]
Note that for $\eta\in A_{\xi},$ we have $|\eta|\sim |\xi-\eta|\sim N$  and consequently $|\theta_1(\xi,\eta)|\sim o(1)$ and $|\theta_3(\xi,\eta)|\sim o(1).$ Therefore, for any $t>0$ fixed we get
\begin{equation}\label{IneqDer}
\left| I_1(\xi,\eta,t)-I_2(\xi,\eta,t)\right|\geq C|t|.
\end{equation}
Also note that, for $|\xi|\sim o(1),$ measure$(A_{\xi})\geq 1.$ Hence, for any fixed $t>0$ and for some positive constant $C_0,$ using (\ref{IneqInt}) and (\ref{IneqDer}), we obtain
\begin{align}\label{ineqc0}
\|I_2((\phi_N,\psi_N),(\phi_N,\psi_N),t)\|^2_{X^s}&\apprge |t|\int_{|\xi|\sim o(1)}(1+|\xi|^2)^s|\varphi(\xi)|^2d\xi \notag\\
& \apprge |t|\int_{|\xi|\sim o(1)} |\xi|^2 d\xi \geq C_0.
\end{align}
By construction $\|(\phi_N,\psi_N)\|_{X^s}\rightarrow 0,$ for any $s<0,$ therefore the last inequality ensures that for any fixed $t>0,$ the application $(u_0,v_0)\mapsto I_2((u_0,v_0),(u_0,v_0),t)$ is not continuous at the origin from $X^s(\mathbb{R})=H^s(\mathbb{R})\times H^{s-1}(\mathbb{R})$ to even $\left(\mathcal{D}'(\mathbb{R})\right)^2.$\\

The rest of the proof follows in a very similar way as in Molinet and Vento \cite{MolVen1, MolVen2} and Panthee \cite{Panthee1}, we make an sketch for the sake of completeness. The idea is to prove that the discontinuity of $(u_0,v_0) \mapsto  I_2((u_0,v_0),(u_0,v_0), t)$ at the origin implies the discontinuity of the flow-map $(u_0,v_0)\mapsto (u(t),v(t)).$ From Theorem \ref{TheoBanquet},  there exist $T > 0$ and $\epsilon > 0$ such that for any $|\epsilon| \leq \epsilon_0$, any $\|(h,h)\|_{L^2(\mathbb{R})\times L^2(\mathbb{R})} \leq1$ and $0\leq t\leq T$, one has
\begin{equation}\label{serieu}
\overrightarrow{u}(\epsilon\overrightarrow{h},t)=\epsilon S(t)\overrightarrow{h}+\sum_{k=2}^{\infty}\epsilon^k I_k({\overrightarrow{h}}^k,t),
\end{equation}
where${ \overrightarrow{h}}^k:= ((h,h), (h,h), \cdots, (h,h)),$ $ \overrightarrow{h}^k\mapsto I_k({\overrightarrow{h}}^k,t)$ is a $k$-linear continuous map from $X^0(\mathbb{R})$ into $C([0,T];X^0(\mathbb{R}))$ and the series converges absolutely in $C([0,T];X^0(\mathbb{R}))$ (See Theorem 3 in Bejenaru and Tao \cite{BejTao1}).\\

From (\ref{serieu}), we have that 
\begin{equation}\label{ineqser1}
u(\epsilon(\phi_N,\psi_N),t)-\epsilon^2I_2((\phi_N,\psi_N), (\phi_N,\psi_N),t)=\epsilon S(t)\overrightarrow{h}+\sum_{k=3}^{\infty}\epsilon^k I_k({\overrightarrow{h}}^k,t),
\end{equation}
Note also that 
\begin{equation}\label{ineqser2}
\|S(t)(\phi_N,\psi_N)\|_{X^s}\leq \|(\phi_N,\psi_N)\|_{X^s}\sim N^s
\end{equation}
and
\begin{equation}\label{ineqser3}
\left\|\sum_{k=3}^{\infty}\epsilon^k I_k({(\phi_N,\psi_N)}^k,t)\right\|_{X^0}\leq \left(\frac{\epsilon_{\ }}{\epsilon_0}\right)^3\sum_{k=3}^{\infty}\epsilon_0^k\| I_k({(\phi_N,\psi_N)}^k,t)\|_{X^0}\leq C\epsilon^3.
\end{equation}
Therefore, from (\ref{ineqser1}), (\ref{ineqser2}) and (\ref{ineqser3}) we get, for any $s<0,$
\[
\sup_{t\in[0,T]}\|u(\epsilon(\phi_N,\psi_N),t)-\epsilon^2I_2((\phi_N,\psi_N), (\phi_N,\psi_N),t)\|_{X^s}\leq O(N^s)+C\epsilon^3.
\]
Now, if we fix $0 < t < 1,$  take $\epsilon$ small enough and then $N$ large enough, and taking into account  (\ref{ineqc0}); the estimate (\ref{ineqser3}) yields that $\epsilon^2I_2((\phi_N,\psi_N), (\phi_N,\psi_N),t)$ is a good approximation of $u(\epsilon(\phi_N,\psi_N),t)$ in $X^s(\mathbb{R})$ for any $s < 0.$\\

If we choose $\epsilon \ll 1 $, from (\ref{ineqser1}), (\ref{ineqser2}) and (\ref{ineqser3}), we get
\begin{align*}
\| u(\epsilon(\phi_N,\psi_N),t)\|_{X^s}&\geq \epsilon^2\|I_2((\phi_N,\psi_N), (\phi_N,\psi_N),t)\|_{X^s}-\epsilon \|S(t)(\phi_N,\psi_N)\|_{X^s}\\
& -\sum_{k=3}^{\infty}\epsilon^k\| I_k({(\phi_N,\psi_N)}^k,t)\|_{X^s}\\
&\geq C_0\epsilon^2-C_1\epsilon^3-C\epsilon N^s\\
&\geq \frac{C_0}{2}\epsilon^2-C\epsilon N^s.
\end{align*}
If we fix the $\epsilon \ll 1$  chosen earlier and choose $N$ large enough, then for any $s < 0,$ the last estimate  yields,
\[\| u(\epsilon(\phi_N,\psi_N),t)\|_{X^s}\geq \frac{C_0}{4}\epsilon^2.\]
Note that, $(u(0,t),v(0,t)\equiv 0$ and $\|(\phi_N,\psi_N)\|_{X^s(\mathbb{R})}\rightarrow 0$ for any $s < 0.$ Therefore, taking  $N\rightarrow \infty$ we conclude that the flow-map $(u_0,v_0)\mapsto(u(t),v(t))$ is discontinuous at the origin from $X^s(\mathbb{R})$ to $C\left([0,1];X^s(\mathbb{R})\right),$ for $s < 0.$ Moreover, as $(\phi_N,\psi_N)\rightharpoonup 0$  in $X^0(\mathbb{R}),$ we also have that the flow-map is discontinuous from $X^0(\mathbb{R}),$ equipped with its weak topology inducted by $X^s(\mathbb{R})$ with values even in $(\mathcal{D}'(\mathbb{R}))^2.$
\endproof

%%%%%%%%%%%%%%%%%%%%% REMARK %%%%%%%%%%%%%%%%%%%%%%%%%%%%%
\begin{remark}\label{remark1}
In the periodic case, i.e., for $x\in\mathbb{T},$ there is analytical well-posedness result for data given in Sobolev spaces without zero Fourier mode ( i.e., with zero $x-$mean) $H^s(\mathbb{T})\times H^{s-1}(\mathbb{T})$,  $s\geq 0$, see Banquet \cite{Banquet1}. Now, for $N \gg 1,$  define the sequence of functions $a_n$ and $b_n$, by
\[a_n=\left \{ 
\begin{aligned}
&1, \ \ |n|\sim N\\
&0, \ \ \text{otherwise}\\
\end{aligned} \right.  \ \ \ \ \  \ \ \ b_n=\left \{ 
\begin{aligned}
& \sqrt{1+N^2}, \ \ |n|\sim N\\
&0, \ \ \text{otherwise}\\
\end{aligned} \right. 
\]
Consider $\widehat{\phi_N},$  $\widehat{\psi_N}$ given by $\widehat{\phi_N}(n)=a_n,$  $\widehat{\psi_N}(n)=b_n,$ then clearly   $\|(\phi_N,\psi_N)\|_{L^2(\mathbb{T})\times H^{-1}(\mathbb{T})}\sim 1,$  and $\|(\phi_N,\psi_N)\|_{H^s(\mathbb{T})\times H^{s-1}(\mathbb{T})}\rightarrow 0,$  for any $s < 0.$ If we proceed with the calculations exactly as above considering the Sobolev space  $H^s(\mathbb{T})\times H^{s-1}(\mathbb{T})$ without zero Fourier mode, we can obtain a similar ill-posedness result for $s < 0,$ in the periodic case too.
\end{remark}
%%%%%%%%%%%%%%% LONG PERIOD LIMIT  %%%%%%%%%%%%%%%%%%%% LONG PERIOD LIMIT %%%%%%%%%%%%%%%%
\section{Long Period limit}
In this section we study the convergence of the periodic solutions of the SRLW equation to the continuous one, on Sobolev spaces. First we study the continuous initial value problem, namely,
\begin{equation}\label{ivpcont}
\left \{ 
\begin{aligned}
&u_t-u_{xxt}+uu_x-v_x=0,\ \ \ \ &x\in\mathbb{R}, \ \ t>0\\
&v_t-u_x=0, & x\in\mathbb{R}, \ \ t>0 \\
&(u(x,0),v(x,0))=(\psi(x),\phi(x)), &x\in\mathbb{R}.
\end{aligned} \right.
\end{equation}
Rewrite the first equation in (\ref{ivpcont}) as $u_t-u_{xxt}=v_x-uu_x,$ then, following Benjamin, Bona and Mahony \cite{benBonMah}, formally solve for $u_t$ to obtain
\[u_t(x,t)=\frac 1{2}\int_{\mathbb{R}}e^{-|x-y|}\left[v_y(y,t)-u(y,t)u_y(y,t)\right]dy.\]
Integrating by parts on the right-hand side of the last equation yields 
\[u_t(x,t)=-\frac 1{2}\int_{\mathbb{R}}\text{sgn}(x-y)e^{-|x-y|}\left[v(y,t)-\tfrac{1}{2}u^2(y,t)\right]dy.\]
Formally, integrating with respect to the temporal variable over $[0,t],$ the last equation and the second equation in (\ref{ivpcont}) one obtains the integral equations
\begin{equation}\label{intequas}
\left \{ 
\begin{aligned}
&u(x,t)=\psi(x)+\int_0^t\int_{\mathbb{R}}K(x-y)\left[v(y,t)-\tfrac{1}{2}u^2(y,t)\right]dy,\\
&v(x,t)=\phi(x)+\int_0^tu_x(x,s)ds, \\
\end{aligned} \right.
\end{equation}
where $K(x)=-\frac{1}{2}\ \text{sgn}(x)e^{-|x|}.$\\

As was pointed out before, Banquet in \cite{Banquet1} established the next result.
\begin{theo}
The initial value problem (\ref{ivpcont}) is locally and globally well-posed in $H^s(\mathbb{R})\times H^{s-1}(\mathbb{R})$ if $s\geq 0.$
\end{theo}
For all the discussion about convergence, we assume that the Fourier transform of a function $f$ is given by 
\[\widehat{f}(\xi)=\int_{\mathbb{R}}f(x)e^{-2\pi i x\xi}dx\]
and the norm of the continuous Sobolev space is
\[\|f\|_{H^s}=\left(\int_{\mathbb{R}}(1+|2\pi \xi|^{2s})|\widehat{f}(\xi)|^2d\xi\right)^{\frac 1{2}}.\]
The following theorem gives spatial decay estimates of solutions $u$ of the Cauchy problem (\ref{ivpcont}) corresponding to similar conditions on the initial data $\psi.$
\begin{theo}
Suppose that $\psi\in H^1(\mathbb{R})$ and $r(x)\psi(x)$ is uniformly bounded, where $r(x)=(1+x^2)^{\sigma}$ with $\sigma>0$  a constant, or $r(x)=e^{\lambda x}$ with $\lambda\in(0,1).$ Then for any $x\in\mathbb{R}$ and $T>0,$ the solution $u$ satisfies 
\begin{equation}\label{desimpcont}
|r(x)u(x,t)|\leq (\|r\psi\|_{\infty}+\kappa T\|r\phi\|_{\infty})e^{ \kappa\left(T+\|(\psi,\phi)\|_{X^{1,2}}\right)t}, \ \ \ \text{for all} \ \ t\in[0,T].
\end{equation}
Here
\[\kappa=\sup_{x\in\mathbb{R}}\int_{\mathbb{R}}\frac{r(x)}{r(y)}e^{-|x-y|}dy.\]
\end{theo}
\proof Define $U(x,t)=r(x)u(x,t)$ and $V(x,t)=r(x)v(x,t).$ Then, the first integral equation in (\ref{intequas}) is equivalent to 
\begin{equation}\label{estu}
U(x,t)=r(x)\psi(x)+r(x)\int_0^t\int_{\mathbb{R}}\frac{K(x-y)}{r(y)}\left[V(y,s)-\frac{1}{2}u(y,s)U(y,s)\right]dyds.
\end{equation}
Using the second equation in (\ref{intequas}) and integrating by parts,  we have for $s\in[0,T]$ and $x\in\mathbb{R}$ that
\[\begin{aligned}
r(x)\int_{\mathbb{R}}\frac{K(x-y)}{r(y)}V(y,s)dy&=\int_{\mathbb{R}}r(x)K(x-y)v(y,s)dy\\
&=\int_{\mathbb{R}}r(x)K(x-y)\left(\phi(y)+\int_0^su_y(y,\tau)d\tau\right)dy\\
&=\int_{\mathbb{R}}r(x)K(x-y)\phi(y)dy+\frac{1}{2}\int_0^s\int_{\mathbb{R}}r(x)e^{-|x-y|}u(y,\tau) dyd\tau\\
&=\int_{\mathbb{R}}\frac{r(x)}{r(y)}K(x-y)r(y)\phi(y)dy+\frac{1}{2}\int_0^s\int_{\mathbb{R}}\frac{r(x)}{r(y)}e^{-|x-y|}U(y,\tau)dyd\tau.
\end{aligned}\] 
Since $s\leq t \leq T,$ we get that
\begin{align}\label{estv}
\left | r(x)\int_0^t\int_{\mathbb{R}}\frac{K(x-y)}{r(y)}V(y,s)dyds\right |&\leq \kappa \int_0^t\left(\|r \phi \|_{\infty}+\int_0^s\|U(\tau)\|_{\infty}d\tau\right)ds\notag \\
&\leq T\kappa\left(\|r \phi \|_{\infty}+\int_0^t\|U(\tau)\|_{\infty}d\tau\right). 
\end{align}
Now, using (\ref{estu}), (\ref{estv}) and elementary considerations we obtain
\[
|U(x,t)|\leq \|r\psi\|_{\infty}+\kappa T\|r\phi\|_{\infty}+\kappa \int_0^t\left(T+\frac{1}{4}\|u(s)\|_{\infty}\right)\|U(s)\|_{\infty}ds.
\]
Since
\[\|u(s)\|_{\infty}\leq \frac{\sqrt{2}}{2}\|u(s)\|_{H^1}\leq \frac{\sqrt{2}}{2}\|(u(s),v(s))\|_{X^{1,2}}=\frac{\sqrt{2}}{2}\|(\psi,\phi)\|_{X^{1,2}},\]
we have that
\[\|U(t)\|_{\infty}\leq \|r\psi\|_{\infty}+\kappa T\|r\phi\|_{\infty}+\kappa \left(T+\|(\psi,\phi)\|_{X^{1,2}}\right)\int_0^t\|U(s)\|_{\infty}ds.\]
Applying Gronwall inequality gives
\[\|U(t)\|_{\infty}\leq \left(\|r\psi\|_{\infty}+\kappa T\|r\phi\|_{\infty}\right)e^{\kappa \left(T+\|(\psi,\phi)\|_{X^{1,2}}\right)t}\]
for all $t\in[0,T]$ and the estimate (\ref{desimpcont}) follows readily. 
\endproof

 Next, we review the problem (\ref{ivpcont}) with periodic initial data. Banquet in \cite{Banquet1} proved the following result.
 \begin{theo}
 The initial value problem (\ref{ivpcont}) is locally and globally well-posed in $H_{per}^s\times H_{per}^{s-1}$ if $s\geq 0.$
 \end{theo}
 The following theorem gives the estimative of the periodic solution $u$ of (\ref{ivpcont}) at $x=\pm l.$ 
 \begin{theo}
 In Cauchy  problem (\ref{ivpcont}), assume that the periodic initial data $(\psi,\phi)\in H^1_l\times L^2_l$ and let $T>0,$
 \[D=l\sum_{n\in\mathbb{Z}}|\psi_n-\psi_{n-1}|, \ \ \ E=l\sum_{n\in\mathbb{Z}}|\phi_n-\phi_{n-1}|,\]
 where $f_n$ is the Fourier coefficients of $f.$ Denote $A=T+\beta(l)\|(\psi,\phi)\|_{X_l^{1,2}}$
  and 
  \[B=\pi T\beta(l)\|(\psi,\phi)\|_{X_l^{1,2}}\left[T+\beta(l)\|(\psi,\phi)\|_{X_l^{1,2}},\right]\]
  with 
  \[\beta(l)=\left(\frac 1{2l}\sum_{n\in\mathbb{Z}}\left( 1+\left|\frac{n\pi}{l}\right|^2\right)^{-1}\right)^{\frac 1{2}}.\]
  Then at $x=\pm l,$ the solution $u(x,t)=\sum_{n\in\mathbb{Z}}u_n(t)e^{\frac{in\pi}{l}x}$ of (\ref{ivpcont}) has the bound 
  \[l|u(l,t)|=l|u(-l,t)|\leq\frac{D+TE}{2}e^{At}+\frac{B}{2A}\left(e^{At}-1\right), \ \ \text{for all}\ \ t\in[0,T].\]
\end{theo}
\proof Take the Fourier transform on the equations in (\ref{ivpcont}), solve for  $u_t^n$ and $v_t^n,$  integrate from $0$ to $t$ and use the initial conditions to obtain 
\[
\left \{ 
\begin{aligned}
&u_n(t)=\psi_n-\int_0^t\frac{\frac{in\pi}{l}}{1+\left|\frac{n\pi}{l}\right|^2}\left[\frac1{2}\sum_{k\in\mathbb{Z}}u_{n-k}(\tau)u_k(\tau)-v_n(\tau)\right]d\tau,\\
&v_n(t)=\phi_n+\int_0^t\frac{in\pi}{l}u_n(\tau)d\tau.\\
\end{aligned} \right.
\]
Let $q_n(t)=l[u_n(t)-u_{n-1}(t)]$ and  $p_n(t)=l[v_n(t)-v_{n-1}(t)],$ for all positive integer $n.$ Then $q_n$ and $p_n$ satisfy the next system of equations
\begin{equation}\label{sistperio}
\left \{ 
\begin{aligned}
q_n(t)&=l[\psi_n-\psi_{n-1}]-\frac{\frac{in\pi}{l}}{1+\left|\frac{n\pi}{l}\right|^2}\int_0^t\left[\frac1{2}\sum_{k\in\mathbb{Z}}u_{n-k}(\tau)q_k(\tau)-p_n(\tau)\right]d\tau\\
&-\frac{i\pi\left[1-\frac{n(n-1)\pi^2}{l^2}\right]}{\left[1+\left|\frac{(n-1)\pi}{l}\right|^2\right]\left[1+\left|\frac{n\pi}{l}\right|^2\right]}\int_0^t\left[\frac1{2}\sum_{k\in\mathbb{Z}}u_{n-1-k}(\tau)u_k(\tau)-v_{n-1}(\tau)\right]d\tau\\
p_n(t)&=l[\phi_n-\phi_{n-1}]+\frac{in\pi}{l}\int_0^tq_n(\tau)d\tau+i\pi \int_0^t u_{n-1}(\tau)d\tau.\\
\end{aligned} \right.
\end{equation}
Define 
\[Q(t)=\sum_{n\in\mathbb{Z}}|q_n(t)| \ \ \ \text{and}\ \ \ P(t)=\sum_{n\in\mathbb{Z}}|p_n(t)|.\]
Using the second equation in (\ref{sistperio}), we get that
 \[
\left|\sum_{n\in\mathbb{Z}}\frac{\frac{in\pi}{l}}{1+\left|\frac{n\pi}{l}\right|^2}\int_0^t p_n(\tau)d\tau\right| \leq TE+T\int_0^tQ(\tau)d\tau+\pi T\int_0^t\sum_{n\in\mathbb{Z}}|u_n(\tau)|d\tau.
 \]
Note that
\[\sum_{n\in\mathbb{Z}}|u_n(\tau)|\leq \beta(l)\|u(t)\|_{H^1_l}\leq \beta(l)\|(u(t),v(t))\|_{X^{1,2}_l}=\beta(l)\|(\psi,\phi)\|_{X^{1,2}_l}\] 
and consequently 
 \begin{equation}\label{ineq1}
 \left|\sum_{n\in\mathbb{Z}}\frac{\frac{in\pi}{l}}{1+\left|\frac{n\pi}{l}\right|^2}\int_0^t p_n(\tau)d\tau\right| \leq TE+T\int_0^tQ(\tau)d\tau +\pi T^2\beta(l)\|(\psi,\phi)\|_{X^{1,2}_l} .
 \end{equation}
 On the other hand, using the fact that $\frac{|1-xy|}{(1+x^2)(1+y^2)}\leq \frac{4}{1+x},$ for all $x,y>0,$ we arrived at
 \begin{align}\label{ineq2}
 \left|\sum_{n\in\mathbb{Z}}\frac{i\pi\left[1-\frac{n(n-1)\pi^2}{l^2}\right]}{\left[1+\left|\frac{(n-1)\pi}{l}\right|^2\right]\left[1+\left|\frac{n\pi}{l}\right|^2\right]}\int_0^tv_{n-1}(\tau)d\tau\right|&\leq 4\pi\int_0^t\sum_{n\in\mathbb{Z}}\left(1+\left|\frac{n\pi}{l}\right|\right)^{-1}|v_n(\tau)|d\tau\notag \\
 &\leq 4\pi\int_0^t\beta(l)\|v(\tau)\|_{L_l^2}d\tau\notag \\
 &\leq 4\pi\beta(l)\int_0^t\|(u(\tau),v(\tau))\|_{X_l^{1,2}}d\tau\notag\\
 &= 4\pi T\beta(l)\|(\psi,\phi)\|_{X_l^{1,2}}.
 \end{align}
 Now, from the first integral equation in (\ref{sistperio}), (\ref{ineq1}) and (\ref{ineq2}), we obtain
 \[
 \begin{aligned}
Q(t)&\leq \int_0^t \left[\frac{1}{2}\sum_{n\in\mathbb{Z}}|u_n(\tau)|Q(\tau)+TQ(\tau)\right]d\tau +D +TE+ \pi T^2\beta(l)\|(\psi,\phi)\|_{X^{1,2}_l}
\\ &+\frac{\pi}{2}\int_0^t \left(\sum_{n\in\mathbb{Z}}|u_n(\tau)|\right)^2d\tau\leq \int_0^t\left[\beta(l)\|(\psi,\phi)\|_{X^{1,2}_l}+T\right]Q(\tau)d\tau+D+TE\\
&+\pi T\beta(l)\|(\psi,\phi)\|_{X^{1,2}_l}\left[T+\beta(l)\|(\psi,\phi)\|_{X^{1,2}_l}\right].
 \end{aligned}
 \]
Therefore
\[Q(t)\leq A\int_0^tQ(\tau)d\tau+B+D+TE.\]
Thus, using the Gronwall inequality, it follows that
\[Q(t)\leq (D+TE)e^{At}+\frac{B}{A}\left(e^{At}-1\right),\ \ \text{for all} \ \ t\in [0,T].\]
 Now, we consider the values of the solution $u$ at $x=\pm l,$ it is seen that
 \[u(l,t)=u(-l,t)=\sum_{n\in\mathbb{Z}}(-1)^nu_n(t)=-\sum_{n\in\mathbb{Z}}\left[u_{2n+1}(t)-u_{2n}(t)\right]=\sum_{n\in\mathbb{Z}}\left[u_{2n}(t)-u_{2n-1}(t)\right]\]
 and hence
 \[l|u(l,t)|=lu(-l,t)|\leq\frac{1}{2}Q(t)\leq \frac{D+TE}{2}e^{At}+\frac{B}{2A}\left(e^{At}-1\right),\]
for all $t\in[0,T],$ which finishes the proof of the theorem.
 \endproof
 \begin{coro}
 In the last theorem, if
 \[D=D(l)=l\sum_{n\in\mathbb{Z}}|\psi_n-\psi_{n-1}|, \ \ \ E=E(l)=l\sum_{n\in\mathbb{Z}}|\phi_n-\phi_{n-1}|,\ \ \ \|\psi\|_{H^1_l} \ \ \ \text{and} \ \ \ \|\phi\|_{L^2_l}\]
  are uniformly bounded as $l\to \infty,$ then on any compact interval $[0,T],$ there exists a constant $C$ independent of $l$ such that
 \[|u(l,t)|=|u(-l,t)|\leq \frac{C}{l}, \ \ \ \text{for all}\ \ \ t\in[0,T].\]
 \end{coro}
 %%%%%%%%%%%%%%% LIMITING RESULTS %%%%%%%%%%%%%%%%%%%%%%%%%%%%%%%%%%%%%%%%%%%%%%%%%%%%%%%%%%%%
 \subsection{Limiting Results}
 For reader's convenience we repeat the initial value problem associated with the SRLW equation, namely
 \begin{equation}\label{ivpcontrep}
\left \{ 
\begin{aligned}
&u_t-u_{xxt}+uu_x-v_x=0,\ \ \ \ &x\in\mathbb{R}, \ \ t>0\\
&v_t-u_x=0, & x\in\mathbb{R}, \ \ t>0 \\
&(u(x,0),v(x,0))=(\psi(x),\phi(x)), &x\in\mathbb{R}.
\end{aligned} \right.
\end{equation}
Suppose that the initial data $\psi$ and $\phi$ are sufficiently nice (See Theorem \ref{printheo} below), so that their Fourier $\widehat{\psi}$ and $\widehat{\phi}$ exist and are continuous. Since we are interested in real functions we have that $\widehat{f}(-\xi)=\overline{\widehat{f}(\xi)}.$ Assuming that $\widehat{f}$ is continuous in $\mathbb{R},$ for $l>0$ we introduced a transform $\mathcal{P}_l$ as follows
\[\mathcal{P}_l(f)(x)=\sum_{n\in\mathbb{Z}}\frac{1}{2l}\widehat{f}\left(\frac{n}{2l}\right)e^{\frac{in\pi}{l}x}.\]      
The function $\mathcal{P}_l(f)$ is formally a real periodic function of period $2l$ because $\widehat{f}(\frac{n}{2l})$ and  $\widehat{f}(\frac{-n}{2l})$ are complex conjugates for any $n.$  Denote  $u_l=u^l$ and  consider the periodic problem
 \begin{equation}\label{ivpperio}
\left \{ 
\begin{aligned}
&u^l_t-u^l_{xxt}+u^lu^l_x-v^l_x=0,\ \ \ \ &x\in\mathbb{R}, \ \ t>0\\
&v^l_t-u^l_x=0, & x\in\mathbb{R}, \ \ t>0 \\
&(u_l(x,0),v_l(x,0))=(\mathcal{P}_l(\psi)(x),\mathcal{P}_l(\phi)(x)), &x\in\mathbb{R}.
\end{aligned} \right.
\end{equation}
The  next proposition was proved by Chen  (see Proposition 3.1 in \cite{ChenHong}),  for $m$ a positive integer, we only note that the argument used to obtain this result also works on the case $m=0.$
 \begin{prop}
 Let $m\in\mathbb{N} \cup \{0\}.$ If $f\in H^m(\mathbb{R})$ and $f,f^{m},\widehat{f^{m}}\in L^1(\mathbb{R}),$ then for any $l>0$ the periodic function $\mathcal{P}_l(f)\in H^m_l,$ and for any $\epsilon >0,$ there exist an $l_{\epsilon}>0$ sufficiently large such that when both $\mathcal{P}_l(f)$ and $f$ are restricted to the interval $(-l,l),$
 \[\|\mathcal{P}_l(f)-f\|_{W_l^{m,2}}<\epsilon, \ \ \ \text{for all} \ \ \ l\geq l_{\epsilon}.\]  
 \end{prop}
 Next we present our principal result on the comparisons between the solutions $(u,v)$ of (\ref{ivpcontrep}) and $(u_l,v_l)$ of (\ref{ivpperio}).
 %%%%%%%%%%%% PRINCIPAL RESULT %%%%%%%%%%%%%%%%%%%%%%%%%%%%%%%%%%%%%%%%%%%%%%%%%%%%%%%
 \begin{theo}\label{printheo}
 Consider the Cauchy problem (\ref{ivpcontrep}), assume that the initial data $(\psi,\phi)$ satisfies: $(\psi,\phi)\in H^1(\mathbb{R})\times L^2(\mathbb{R}),$ $\psi',  x\phi(x)\in L^1(\mathbb{R})$ and there exists $s>1$ such that
 \[(1+x^2)^{\frac{s}{2}}\psi(x)\in L^{\infty}(\mathbb{R}).\]
 Moreover, assume that the Fourier transform of $\psi$ and $\phi$ satisfy: 
 \[(1+\xi^2)^{\frac{1}{2}}\widehat{\psi}(\xi), \widehat{\phi}, \widehat{x\psi(x)}, \widehat{x\phi(x)}\in L^1(\mathbb{R}).\]
 Then, the solutions $u$ of (\ref{ivpcontrep}) and the solution $u_l$ of (\ref{ivpperio}), when restricted to the spatial interval $(-l,l),$ have the relation
 \[\lim_{l\to\infty}\|(u_l(t),v_l(t))-(u(t),v(t))\|_{X^{1,2}_l}=0.\]
 This convergence in uniformly on any compact interval $[0,T].$ 
  \end{theo}
  \proof 
 Introduce new dependent variables $z, \eta$ and $w$ as follows
 \[z(x,t)=u_l(x,t)-u(x,t), \ \ \ \eta(x,t)=v_l(x,t)-v(x,t)\ \ \text{and}\ \ w(x,t)=z(x,t)+\varphi(x,t),\]
 where $\varphi(x,t)=-z(l,t)\varphi_+(x)-z(-l,t)\varphi_-(x).$
 Here 
 \[\varphi_-(x)=\frac{e^{l-x}-e^{-l+x}}{e^{2l}-e^{-2l}}\ \ \ \text{and}\ \ \ \varphi_+(x)=\varphi_-(-x).\]
 Then, $(w(t),\eta(t))\in X^{1,2}(-l,l)$ and satisfy the initial value boundary problem
\begin{equation}\label{ivpmod}
\left \{ 
\begin{aligned}
&w_t-w_{xxt}+ww_x-\eta_x=\left(\varphi u-\frac{1}{2}\varphi^2-(u-\varphi)w\right)_x,\ \ \ \ & x\in [-l,l], \ \ t>0\\
&\eta_t-w_x=-\varphi_x, & x\in[-l,l], \ \ t>0 \\
&w(x,0)=:w_0(x)=\mathcal{P}_l(\psi)(x)-\psi(x)-\left[\mathcal{P}_l(\psi)(l)-\psi(l)\right]\varphi_+(x)\\
&\hspace{2.4cm} -\left[\mathcal{P}_l(\psi)(-l)-\psi(-l)\right]\varphi_-(x), & x\in[-l,l]\\
&w(-l,t)=w(l,t)=0, & t\geq 0\\
&\eta(x,0)=:\eta_0(x)=\mathcal{P}_l(\psi)(x)-\psi(x), & x\in[-l,l].
\end{aligned} \right.
\end{equation}
Multiply the first equation in (\ref{ivpmod}) by $2w$ and integrate over $[-l,l]$ with respect to $x;$ after integration by parts there appears
\begin{align}\label{ineqfin}
\frac{d}{dt}\|(w,\eta)(t)\|^2_{X^{1,2}_l}&=\frac{d}{dt}\int_{-l}^l\left[w^2(x,t)+w^2_x(x,t)+\eta^2(x,t)\right]dx\notag\\
&=2\int_{-l}^l(u-\varphi)ww_xdx-\int_0^t\left(2\varphi u-\varphi^2\right)w_xdx-\int_0^t\varphi_x\eta dx\notag\\
&\leq \|u(t)-\varphi(t)\|_{L^{\infty}_l}\|w(t)\|^2_{W^{1,2}_l}+\|2\varphi(t)u(t)-\varphi^2(t)\|_{L^2_l}\|w_x(t)\|_{L^2_l}\notag\\
&\ \ \ \ +\|\varphi_x(t)\|_{L^2_l}\|\eta(t)\|_{L^2_l}\notag\\
&\leq \left[\|\varphi(t)\|_{L^2_l}\|2u(t)-\varphi(t)\|_{L^{\infty}_l}+\|\varphi_x(t)\|_{L^2_l}\right]\|(w,\eta)(t)\|_{X^{1,2}_l}\\ 
&\ \ \ \ +\|u(t)-\varphi(t)\|_{L^{\infty}_l}\|(w,\eta)(t)\|^2_{X^{1,2}_l}.\notag
\end{align}
 Following the ideas of Chen in \cite{ChenHong}, we have that $\|u(t)-\varphi(t)\|_{L^{\infty}_l}\leq 2\|u(t)\|_{\infty}+|u_l(l,t)|.$  But
 \[\|u(t)\|_{\infty}\leq \frac{1}{\sqrt{2}}\|u(t)\|_{H^1}\leq \frac{1}{\sqrt{2}}\|(u,v)(t)\|_{X^{1,2}}= \frac{1}{\sqrt{2}}\|(\psi,\phi)\|_{X^{1,2}}\]
 and 
 \[|u_l(x,t)|\leq \beta(l)\|u_l(t)\|_{H^1_l}\leq \beta(l)\|(u_l,v_l)(t)\|_{X^{1,2}_l}=\beta(l)\|(\mathcal{P}_l(\psi),\mathcal{P}_l(\phi))\|_{X^{1,2}_l}.\]
 Therefore
 \begin{equation}\label{ineqfin1}
 \|u(t)-\varphi(t)\|_{L^{\infty}_l}\leq \sqrt{2}\|(\psi,\phi)\|_{X^{1,2}}+\beta(l)\|(\mathcal{P}_l(\psi),\mathcal{P}_l(\phi))\|_{X^{1,2}_l}.
 \end{equation}
 It is not difficult to show that
 \[\|2u(t)-\varphi(t)\|_{L^{\infty}_l}\leq 2\sqrt{2}\|(\psi,\phi)\|_{X^{1,2}}+2\beta(l)\|(\mathcal{P}_l(\psi),\mathcal{P}_l(\phi))\|_{X^{1,2}_l}.\]
 Furthermore
 \begin{equation}\label{ineqfin2}
 \|\varphi(t)\|_{L^2_l}\leq \frac{3}{2}U(l,t),
 \end{equation}
  where $U(l,t)=|u(l,t)|+|u(-l,t)|+|u_l(l,t)|.$ Since 
 \[\|\varphi'_++\varphi_-'\|^2_{L^2_l}=\left(e^{2l}-e^{-2l}\right)^2 \left[8l-4l\left(e^{2l}+e^{-2l}\right)+2\left(e^{-2l}-e^{2l}\right)+e^{4l}-e^{-4l}\right]<1,\]
 for all $l>0,$ we have that 
 \begin{equation}\label{ineqfin3}
 \|\varphi_x(t)\|_{L^2_l}\leq \left[\max\{|u(l,t)|,|u(-l,t)|\}+|u_l(l,t)|\right]\|\varphi'_++\varphi_-'\|^2_{L^2_l}\leq U(l,t),
\end{equation} 
Denote $\gamma=\gamma(l)=\frac{1}{2}\left(\sqrt{2}\|(\psi,\phi)\|_{X^{1,2}}+2\beta(l)\|(\mathcal{P}_l(\psi),\mathcal{P}_l(\phi))\|_{X^{1,2}_l}\right).$ From (\ref{ineqfin1}), (\ref{ineqfin2}) and (\ref{ineqfin3}), we can rewrite (\ref{ineqfin}) as
 \[\frac{d}{dt}\|(w,\eta)(t)\|^2_{X^{1,2}_l}\leq 2\gamma \|(w,\eta)(t)\|^2_{X^{1,2}_l}+U(l,t)(6\gamma +1)\|(w,\eta)(t)\|_{X^{1,2}_l}.\]
Hence, Gronwall inequality provides the following estimative
 \[
 \|(w,\eta)(t)\|_{X^{1,2}_l}\leq \|(w_0,\eta_0)\|_{X^{1,2}_l}e^{\gamma t} +\frac{6\gamma-1}{2}\int_0^tU(l,s)e^{\gamma(t-s)}ds
 .\]
 Now, from the definition of $w$ and $\eta,$ we obtain
 \[\|(u_l(t),v_l(t))-(u(t),v(t))\|_{X^{1,2}_l}\leq \|\varphi(t)\|_{W^{1,2}_l}+\sqrt{2} \|(w,\eta)(t)\|_{X^{1,2}_l}.\]
 Since $\|\varphi(t)\|_{W^{1,2}_l}\leq 2(|z(l,t)|+|z(-l,t)|),$ we arrived at 
 \[
\begin{aligned}
\|(u_l(t),v_l(t))-&(u(t),v(t))\|_{X^{1,2}_l}\leq 2 (|z(l,t)|+|z(-l,t)|)+\frac{6\gamma-1}{2}\int_0^tU(l,s)e^{\gamma(t-s)}ds \\
&+\sqrt{2}\left(\|\mathcal{P}_l(\psi)-\psi)\|_{W^{1,2}_l}+ 2(|z(l,t)|+|z(-l,t)|) +\|\mathcal{P}_l(\phi)-\phi)\|_{L^2_l} \right)e^{\gamma t}
\end{aligned}
 \]
 The rest of the proof  follows in the same way as Theorem 3.3 in Chen \cite{ChenHong}, we only prove that $D(l)$ and $E(l)$ are uniformly bounded as $l\to\infty$. Indeed, from $\widehat{x\psi(x)}\in L^1(\mathbb{R})$ we obtain that $\widehat{\psi}'(\xi)$ lies in $L^1(\mathbb{R}).$ Hence 
 \begin{equation}\label{ineqintfin}
 \int_{\mathbb{R}}|\widehat{\psi}'(\xi)|d\xi=\lim_{l\to\infty}\sum_{n\in\mathbb{Z}}\frac{1}{2l}|\widehat{\psi}'(\xi_n)|
 \end{equation}
 for any $\xi_n\in\left[\frac{n-1}{2l}, \frac{n}{2l}\right].$ Now, consider $D(l)$ given above, by the Mean Value Theorem, for each $n,$ there exists $\theta_n\in[0,1],$ such that 
 \[D(l)=\sum_{n\in\mathbb{Z}}\left|\widehat{\psi}\left(\frac{n}{2l}\right)-\widehat{\psi}\left(\frac{n-1}{2l}\right)\right|=\sum_{n\in\mathbb{Z}}\frac{1}{2l}\left|\widehat{\psi}\left(\frac{n-1}{2l}+\frac{\theta_n}{l}\right)\right|.\]
 It is seen that $D(l)$ is uniformly bounded as $l\to\infty$ from (\ref{ineqintfin}). Similarly we get that $E(l)$ is uniformly bounded as $l\to\infty,$ which finishes the proof of the theorem.
\endproof

\textbf{Acknowledgements:} The author is grateful to Hongqui Chen for  helpful discussion on the subject.
%%%%%%%%%%%%%%%%%%%%%%%references%%%%%%%%%%%%references%%%%%%%%%%%%%%%%%%%%%%%%%%%%%%%%%%%%%%%%

\end{document}